\documentclass[12pt]{article}
\usepackage{amssymb}
\usepackage{amsfonts}
\usepackage{amsthm}
\usepackage{amsmath}
\usepackage[dvips]{graphicx}

\newcommand{\vw}{\begin{equation}}
\newcommand{\ww}{\end{equation}}

\newtheorem{definition}{Definition}[section]
\newtheorem{theorem}{Theorem}[section]
\newtheorem{lemma}[theorem]{Lemma}
\newtheorem{remark}[theorem]{Remark}
\numberwithin{equation}{section}
\renewcommand{\O}{\mathcal{O}}

\newcommand{\rd}{D}

\def \reff#1{(\ref{#1})}

\begin{document}

\title{$q$-Discrete Painlev\'{e} equations for recurrence coefficients of modified $q$-Freud orthogonal polynomials}
\author{Lies Boelen$^*$, Christophe Smet$^*$, Walter Van Assche\footnote{email: \{lies.boelen\}, \{christophe\} and \{walter\}@wis.kuleuven.be}
 \\Katholieke Universiteit Leuven\\ Department of Mathematics
 \\Celestijnenlaan 200B-box 2400\\3001 Leuven, Belgium
} \maketitle

\begin{abstract}\noindent We present an asymmetric $q$-Painlev\'e equation. We will derive
this using $q$-orthogonal polynomials with respect to generalized
Freud weights: their recurrence coefficients will obey this
$q$-Painlev\'e equation (up to a simple transformation). We will
show a stable method of computing a special solution which gives the
recurrence coefficients. We establish a connection between the
newfound equation and $\alpha$-$q$-P$_{\rm V}.$\end{abstract}

\section{Introduction}\label{section:intro}

\noindent The aim of this work is to identify new $q$-discrete
Painlev\'e equations. What exactly are the discrete Painlev\'e
equations (d-P)? One could state  that a d-P is a second-order,
nonautonomous integrable mapping having one of the celebrated
Painlev\'e equations for a continuous limit. This description does
not quite capture the whole story, so let's look at some important
notions about discrete Painlev\'e equations.

\noindent Although some discrete Painlev\'e equations were found
\emph{avant la lettre}, the first steps were made by Br\'ezin and
Kazakov \cite{brezkaz} who found an equation (now known as
d-P$_\textrm{I}$) and computed its continuous limit, the continuous
Painlev\'e I. The definition of {\it singularity confinement}
\cite{sc} was maybe the most important step in the evolution of the
field: it is the discrete analogue of the Painlev\'e property and
describes the behaviour of singularities throughout the evolution of
a discrete equation:
\begin{definition}[\textsc{Singularity Confinement Property}]
Consider a difference equation with independent variable $n$ and
dependent variable $x_n$. If $x_n$ is such that it gives rise to a
singularity for $x_{n+1}$, then there exists a $p \in \mathbb{N}$
such that the singularity is confined to $x_{n+1},\ldots, x_{n+p}$
and  $x_{n+p+1}$ depends only on $x_{n-1}, x_{n-2}, \ldots$
\end{definition}

\noindent Using only this criterion of singularity confinement,
discrete analogues for (continuous) P$_{\rm III}$, P$_{\rm IV}$ and
P$_{\rm V}$ were identified \cite{sc}, \cite{tweede}. It's striking
that in the analogue of P$_\mathrm{III}$ we see the independent
variable $n$ entering the equation only in an exponential way:
$$
x_{n+1}x_{n-1} = \frac{(x_n + \alpha)(x_n+\beta)}{(\gamma q^n
+1)(\delta q^n+1)}
$$
with $\alpha, \beta, \gamma, \delta$ constants. Therefore, it is
considered to be a so-called \emph{$q$-discrete} equation as the
nature of the equation is multiplicative rather than additive like,
e.g., d-P$_\mathrm{I}$
$$
x_{n-1}+x_n+x_{n+1}= \frac{z_n+ \gamma (-1)^n}{x_n} +\delta
$$
with  $z_n = \alpha n+\beta$ and $\alpha, \beta, \gamma, \delta$
constants. This equation shows an odd-even dependence through the
factor $(-1)^n$ when $\gamma \neq 0$. In this case we can introduce
new variables $u_n = x_{2n}$, $v_n=x_{2n+1}$ that lead to the system
$$
\left\{
\begin{array}{rcccl}
v_{n-1} + u_n + v_n &=& \frac{2\alpha n +\beta + \gamma}{u_n}
+\delta  &=& \frac{z_n+\gamma}{u_n}+\delta\\
u_n + v_n+u_{n+1} &=& \frac{2\alpha (n+1) + \beta
-\gamma}{v_n}+\delta  &=& \frac{z_{n+1}-\gamma}{v_n} +\delta
\end{array}\right.
$$
with $z_n = 2n\alpha +\beta$. This system is known as
$\alpha-$d-P$_\textrm{I}$, an \emph{asymmetric} discrete Painlev\'e
equation.  A list of a few important discrete Painlev\'e equations
was compiled by Peter Clarkson, and can be found in \cite{walter}.
See also \cite{gramram} for an overview of discrete Painlev\'e
equations.\vspace{8mm}

\noindent The link between orthonormal polynomials and discrete
Painlev\'e equations is well established (\cite{ifk}, \cite{magnus},
\cite{walter}). Given a positive measure $w$ on a set $A \subset
\mathbb{R}$ and assuming all the moments for $w$ exist, i.e. $|
\int_A x^k w(x) dx| < \infty$, we denote with $\{p_n\}$ the set of
orthonormal polynomials with respect to $w$:
$$
\int_{A} p_n(x) p_m(x) w(x) dx = \delta_{mn}\hspace{10mm} m, n \geq
0.
$$
These orthonormal polynomials are unique if we choose the leading
coefficient to be positive.

\noindent Orthonormal polynomials satisfy a three-term recurrence
relation of the form
$$
xp_n(x) = a_{n+1}p_{n+1}(x) + b_n p_n(x) + a_n
p_{n-1}(x)\hspace{10mm} n\geq 0
$$
with $p_{-1}=0$. Here we can find the recurrence coefficients $a_n$
and $b_n$ as the coefficients showing up in the Fourier expansion of
$xp_n(x)$:
\begin{eqnarray}
\label{eq:intro_an}a_n &=& \int_A x p_n(x) p_{n-1}(x) w(x) dx,\\
\label{eq:intro_bn}b_n &=& \int_A x p_n^2(x) w(x) dx.
\end{eqnarray}

\noindent If we write $p_n(x) = \gamma_n x^n + \ldots$, we can find
$a_n = \frac{\gamma_{n-1}}{\gamma_n}$ by comparing the leading
coefficients of both sides of the recurrence relation. Since we
chose to take positive leading coefficients, $a_n$ will be positive
as well.

\noindent During the last few decades, many  generalizations of
well-known weight functions have been shown to give rise to
discrete Painlev\'e equations for the recurrence coefficients
$a_n$ and $b_n$.

\noindent As far as we can see, the very first appearance of a d-P
was in the context of orthogonal polynomials, when Shohat
\cite{shohat} found a nonlinear recurrence relation for the
recurrence coefficients, which is now known as d-P$_\textrm{I}$.
Fokas, Its and Kitaev found the connection of d-P$_{\rm I}$ to Freud
weights in \cite{ifk}. Nijhoff \cite{nijhoff} found a $q$-discrete
Painlev\'e equation in the context of orthogonal polynomials: he
considered a $q$-generalization of the Hermite polynomials on the
exponential lattice and found non-linear recurrences of order $>$2.
A second order recurrence relation was found by the third author
\cite{walter} which he called q-P$_{\rm I}$. We will study a slight
extension of this Freud weight to recover asymmetric $q$-Painlev\'e
equations.

\noindent In particular we consider in Section 3 the weight
$$
w(x) = \frac{(q^4x^4;q^4)_\infty |x|^\alpha}{(1-q^4)^{\alpha/4}}
\hspace{15mm} \alpha> -1
$$
on the $q$-exponential lattice $L = \{\pm q^n, n \in \mathbb{N}\}$,
with $q \in (0,1)$. We will show that $b_n = 0$ and $y_n =
a_n^2q^{1-n}$ satisfies
\begin{displaymath}
q^{n-\alpha}(y_ny_{n+1}+q^\alpha)(y_ny_{n-1}+q^\alpha) = \left\{
\begin{array}{ll} q^\alpha - q^{-\alpha}y_n^2 & n \textrm{ even}\\
1 - y_n^2 & n \textrm{ odd}
\end{array} \right.
\end{displaymath}
which we will call q-P$_I$. The case $\alpha = 0$ was already
obtained in \cite{walter}. In Section 4 we consider the more general
weight
$$
w(x) = \frac{|x|^\alpha
(q^2x^2;q^2)_\infty(cq^2x^2;q^2)_\infty}{(1-q^4)^{\alpha/4}},
\hspace{10mm} \alpha >-1, c \leq 0.
$$
We show that $y_n=a_n^2q^{1-n}$ now satisfies
\begin{displaymath}
q^{n-\alpha}(-cy_ny_{n+1}+q^\alpha)(-cy_ny_{n-1}+q^\alpha) = \left\{
\begin{array}{ll} (q^\alpha - y_n)(q^\alpha - cy_n)q^{-\alpha} & n \textrm{ even}\\
(1 - y_n)(1-cy_n) & n \textrm{ odd.}
\end{array} \right.
\end{displaymath}
In Section 5 we relate these equations to $\alpha$-$q$-P$_{\rm V}$.
In Section 6 we show a stable method of computing the recurrence
coefficients.

\section{Preliminaries}\label{section:q}
We consider orthogonal polynomials on the exponential lattice
$$L=\{\pm q^n| n \in \mathbb{N}\}, \hspace{10mm} 0<q<1.$$ The orthonormality condition
with respect to a weight $w$ on $L$ is,
\begin{equation}\label{eq:orthonormality}\int_{-1}^1 p_n(x) p_m(x) w(x) d_qx =
\delta_{mn}.\end{equation} The $q$-integral is defined as the sum
$$
\int_{-1}^1 f(x) d_qx = (1-q) \sum_{k=0}^\infty f(q^k)q^k +
(1-q)\sum_{k=0}^\infty f(-q^k)q^k.
$$
The $q$-difference operator $D_q$ will play an important role in our
results:
$$
D_qf(x)= \left\{ \begin{array}{ll}
\displaystyle\frac{f(qx)-f(x)}{x(q-1)} & \textrm{ if } x \neq 0\\
f'(0) & \textrm{ if } x= 0.\end{array}\right.
$$
We will consider even weights $w(x) = w(-x)$. This implies that the
orthonormal polynomials $p_n$ associated with the weight $w$ will
satisfy the symmetry property
\begin{equation}\label{eq:symmetry}p_n(-x) = (-1)^n p_n(x).\end{equation} The
orthonormal polynomials of even degree will be even, those of odd
degree will be odd. The recurrence relation will then take the form
\begin{equation}\label{eq:recrel_algemeen} xp_n(x) = a_{n+1}
p_{n+1}(x) + a_n p_{n-1}(x),\end{equation}
i.e., all coefficients $b_n$ are equal to zero.\\

\noindent We will use the following technique (see, e.g.,
\cite{walter}) to identify $q$-discrete Painlev\'e equations:
\begin{itemize}
\item Given an even weight $w$, let $p_n$ denote the orthonormal
polynomials with respect to $w$.
\item Find the Fourier expansion of the polynomial $D_qp_n$.
\item Compare the coefficients of this polynomial and its Fourier
expansion. This will give rise to a set of equations.
\item After a change of variable, we find a $q$-discrete Painlev\'e
equation. We check the property of singularity confinement, and try
to find out what happens when $q \to 1$.
\end{itemize}

\section{Modified $q$-Freud polynomials}\label{section:nog_geen_c}
We consider the weight
\begin{equation}\label{eq:weight_def}w(x)=\frac{(q^4x^4;q^4)_\infty
|x|^\alpha}{(1-q^4)^{\alpha/4}}\end{equation} on the $q$-exponential
lattice $L$, with $\alpha>-1$. The $q$-Pochhammer symbol is defined
as
$$(a;q)_n=\prod_{j=0}^{n-1}\left(1-aq^j\right),\qquad (a;q)_\infty=\prod_{j=0}^{\infty}\left(1-aq^j\right).$$  Observe that, in terms of
the $q$-exponential function $E_q(z)=(-z;q)_\infty$, we have
\begin{equation*}\label{eq:weight_is_q-exp}w\left(\sqrt[4]{1-q^4}x\right)=|x|^\alpha
E_{q^4}\left(-(1-q^4)q^4x^4\right)\end{equation*} and hence
$w(\sqrt[4]{1-q^4}x)\rightarrow |x|^\alpha e^{-x^4}$ when
$q\rightarrow 1$ (see, e.g., \cite{Lubinsky}).  This limit relation
is the only reason for the presence of the (constant) denominator in
$w$. So this weight can be called a $q$-analogue of the modified
Freud weight $|x|^\alpha e^{-x^4}$. It is easy to check that this
weight satisfies the Pearson equation
\begin{equation}\label{eq:pearson}w\left(\frac{x}{q}\right)=\frac{(1-x^4)}{q^\alpha}w(x).\end{equation}
Let $p_n$ be the orthonormal polynomials associated to $w$ as in
\reff{eq:orthonormality}.  Because of \reff{eq:symmetry} we can put
$p_n(x)=\gamma_nx^n+\delta_nx^{n-2}+...$, and we have the following
\begin{lemma}\label{lemma:gammaquotienten}For orthogonal polynomials
with an even weight one has
\[a_n=\frac{\gamma_{n-1}}{\gamma_n}\qquad and\qquad -\sum_{j=1}^{n-1}a_j^2=\frac{\delta_n}{\gamma_n}.\]\end{lemma}
\begin{proof}These two equations can be obtained by comparing, respectively, the $x^{n+1}$ and $x^{n-1}$ terms in the recurrence relation \reff{eq:recrel_algemeen}.\end{proof}
\begin{lemma}\label{lemma:structure_relation}For even $n$, the polynomials $p_n$  with weight \reff{eq:weight_def} satisfy the following structure
relation:
\begin{equation}\label{eq:structure_even}\rd_qp_n(x)=\frac{B_n}{1-q}p_{n-1}(x)+\frac{A_n}{1-q}p_{n-3}(x)\end{equation}
with
\begin{equation}\label{eq:A_n_even}A_n=\frac{a_na_{n-1}a_{n-2}}{q^{\alpha+n-3}}\end{equation}
\begin{equation}\label{eq:B_n_even}B_n=\frac{a_n}{q^{\alpha+n-1}}\left(\sum_{j=1}^{n+1}a_j^2-q^2\sum_{j=1}^{n-2}a_j^2\right).\end{equation}
For odd $n$, the polynomials $p_n$ satisfy
\begin{equation}\label{eq:structure_odd}\rd_qp_n(x)=\frac{B_n}{1-q}p_{n-1}(x)+\frac{A_n}{1-q}p_{n-3}(x)+ {\rm lower} \: {\rm order} \: {\rm terms}
\end{equation}
with
\begin{equation}\label{eq:A_n_odd}A_n=\frac{a_na_{n-1}a_{n-2}}{q^{\alpha+n-3}}-\left(1-q^{-\alpha}\right)\frac{a_{n-1}}{a_na_{n-2}}\end{equation}
\begin{equation}\label{eq:B_n_odd}B_n=\frac{a_n}{q^{\alpha+n-1}}\left(\sum_{j=1}^{n+1}a_j^2-q^2\sum_{j=1}^{n-2}a_j^2\right)+\frac{1-q^{-\alpha}}{a_n}.\end{equation}
\end{lemma}
\begin{remark}Notice the difference between the closed expression
of the structure relation for even $n$, and the presence of all even
lower degree polynomials for odd $n$.  Despite this difference, the
resulting Painlev\'e equations for even and odd $n$ will have a very
similar structure.\end{remark}
\begin{proof}If we expand $\rd_qp_n$ into a Fourier series, we
obtain \[\rd_qp_n(x)=\sum_{j=0}^{n-1}a_{j,n}p_j(x)\] with
\[a_{j,n}=\int_{-1}^1\rd_qp_n(x)p_j(x)w(x)d_q(x).\]  The symmetry
relation \reff{eq:symmetry} shows that $a_{j,n}=0$ if $n-j$ is even.
For $n-j$ odd we get
\begin{eqnarray*}a_{j,n}&=&2(1-q)\sum_{k=0}^\infty\rd_qp_n(q^k)p_j(q^k)w(q^k)q^k\\
&=&-2\sum_{k=0}^\infty\left(p_n(q^{k+1})-p_n(q^k)\right)p_j(q^k)w(q^k)\\
&=&-2\sum_{k=0}^\infty p_n(q^{k+1})p_j(q^k)w(q^k)+2\sum_{k=0}^\infty
p_n(q^k)p_j(q^k)w(q^k).\end{eqnarray*} Since $n-j$ is odd,
\reff{eq:symmetry} implies that both sums are finite.  Now we
perform a shift in the summation index of the first sum, so that it
contains the expression $w(q^{k-1})$, and we apply the Pearson
equation \reff{eq:pearson} on this.  We can recognize all sums as
$q$-integrals, and we obtain
\begin{eqnarray}\nonumber a_{j,n}&=&-\frac{q^{-\alpha}}{1-q}\int_{-1}^1p_n(x)p_j(x/q)\frac{w(x)}{x}d_qx+\frac{q^{-\alpha}}{1-q}\int_{-1}^1p_n(x)p_j(x/q)x^3w(x)d_qx\\\label{eq:a_jn_drie_integralen}&&+\frac{1}{1-q}\int_{-1}^1p_n(x)p_j(x)\frac{w(x)}{x}d_qx.\end{eqnarray}
From now on we have to make a distinction based on the parity of
$n$.\\
\underline{Case 1: $n$ is even.}  Since $n-j$ is odd, we know from
\reff{eq:symmetry} that $p_j(x)$ and $p_j(x/q)$ are odd polynomials,
hence $\frac{p_j(x)}{x}$ and $\frac{p_j(x/q)}{x}$ are polynomials,
and the orthogonality relations \reff{eq:orthonormality} imply that
the first and the third integral in \reff{eq:a_jn_drie_integralen}
vanish.  To obtain the value of the second integral, we write
$p_j(x/q)x^3$ as a linear combination of the orthonormal polynomials
$p_k(x)$: an easy calculation shows that
\begin{equation*}p_j(x/q)x^3=\frac{\gamma_j}{\gamma_{j+3}q^j}p_{j+3}(x)+\left(\frac{\delta_j}{\gamma_{j+1}q^{j-2}}-\frac{\gamma_j\delta_{j+3}}{\gamma_{j+3}\gamma_{j+1}q^j}\right)p_{j+1}(x)+ {\rm lower} \: {\rm order} \: {\rm terms}\end{equation*}
Using Lemma~\ref{lemma:gammaquotienten} and the orthonormality
relations \reff{eq:orthonormality} we
obtain \reff{eq:structure_even}-\reff{eq:B_n_even}.\\
\underline{Case 2: $n$ is odd.}  The second integral in
\reff{eq:a_jn_drie_integralen} can be computed in exactly the same
way as in the previous case.  However, this time the contributions
of the first and the third integral do not vanish.  For the third
integral, we need to write $p_n(x)/x$ as a linear combination of the
orthonormal polynomials $p_k(x)$.  This can be done because $n$ is
odd, and hence $p_n$ is an odd polynomial.  This yields
\begin{equation*}\frac{p_n(x)}{x}=\frac{\gamma_n}{\gamma_{n-1}}p_{n-1}(x)+\left(\frac{\delta_n}{\gamma_{n-3}}-\frac{\gamma_n\delta_{n-1}}{\gamma_{n-1}\gamma_{n-3}}\right)p_{n-3}(x)+{\rm lower} \: {\rm order} \: {\rm terms}\end{equation*}  As for the first integral in
\reff{eq:a_jn_drie_integralen}, we can use the orthogonality
relations \reff{eq:orthonormality} to write
\begin{eqnarray*}\int_{-1}^1p_n(x)p_j(x/q)\frac{w(x)}{x}d_qx&=&\int_{-1}^1p_n(x)\frac{p_j(x/q)-p_j(x)+p_j(x)}{x}w(x)d_qx\\
&=&\int_{-1}^1p_n(x)p_j(x)\frac{w(x)}{x}d_qx,\end{eqnarray*} which
is the same as the third integral in \reff{eq:a_jn_drie_integralen}.
Combining these results with Lemma~\ref{lemma:gammaquotienten} we
obtain \reff{eq:structure_odd}-\reff{eq:B_n_odd}.\end{proof}Now we
can use these structure relations
\reff{eq:structure_even}-\reff{eq:B_n_odd} to obtain relations
between the recurrence coefficients $a_n$.  Comparing coefficients
of $x^{n-1}$ and $x^{n-3}$ in these structure relations and using
Lemma~\ref{lemma:gammaquotienten} gives
\begin{equation}\label{eq:x^(n-3)_even}a_n^2a_{n-1}^2a_{n-2}^2=q^{n+\alpha-3}\left(q^{n-2}(1-q^2)\sum_{j=1}^{n-2}a_j^2-(1-q^{n-2})a_{n-1}^2\right)\end{equation}
and
\begin{equation}\label{eq:x^(n-1)_even}a_n^2\left(a_{n+1}^2+a_n^2+a_{n-1}^2+(1-q^2)\sum_{j=1}^{n-2}a_j^2\right)=(1-q^n)q^{\alpha+n-1}\end{equation}
for even $n$, and
\begin{equation}\label{eq:x^(n-3)_odd}q^{-\alpha-n+3}a_n^2a_{n-1}^2a_{n-2}^2=q^{n-2}(1-q^2)\sum_{j=1}^{n-2}a_j^2-\left(q^{-\alpha}-q^{n-2}\right)a_{n-1}^2\end{equation}
and
\begin{equation}\label{eq:x^(n-1)_odd}a_n^2\left(a_{n+1}^2+a_n^2+a_{n-1}^2+(1-q^2)\sum_{j=1}^{n-2}a_j^2\right)=\left(q^{-\alpha}-q^n\right)q^{\alpha+n-1}\end{equation}
for odd $n$.  The aim is to obtain a Painlev\'{e}-type recurrence
relation between $a_{n+1}$, $a_n$ and $a_{n-1}$.  Hence to get rid
of the $a_{n-2}$ in \reff{eq:x^(n-3)_even} and
\reff{eq:x^(n-3)_odd}, we replace $n$ by $n+1$.  Keep in mind that
this changes the parity of $n$.  The sum $\sum_{j=1}^{n-1}a_j^2$
that arises in this way, should then be seen as
$\sum_{j=1}^{n-2}a_j^2+a_{n-1}^2$, where we can use
\reff{eq:x^(n-1)_odd} and \reff{eq:x^(n-1)_even} respectively to
write these sums as a function of $a_{n+1}$, $a_n$ and $a_{n-1}$
only.  These manipulations yield the following recurrence
relations:\begin{equation}\label{eq:recrel_a_n_evenodd}\left\{\begin{array}{rl}a_n^2\left(a_{n+1}^2+q^{-\alpha-n+1}a_n^2+q^2a_{n-1}^2+q^{-2n-\alpha+3}a_{n+1}^2a_n^2a_{n-1}^2\right)=\left(1-q^n\right)q^{\alpha+n-1},&n\:
{\rm
even}\\a_n^2\left(a_{n+1}^2+q^{-n+1}a_n^2+q^2a_{n-1}^2+q^{-2n-\alpha+3}a_{n+1}^2a_n^2a_{n-1}^2\right)=\left(q^{-\alpha}-q^n\right)q^{\alpha+n-1},&n\:
{\rm odd.}\end{array}\right.\end{equation}This is a $q$-deformation
of the discrete Painlev\'{e} I equation ${\rm d}$-${\rm P}_{\rm I}$,
which in its most general form is given by
\begin{equation}\label{eq:dPI}x_{n+1}+x_n+x_{n-1}=\frac{an+b+c(-1)^n}{x_n}+d.\end{equation}
Indeed, if we take $x_n=a_n^2/\sqrt{1-q^4}$ and we let $q$ tend to
1, we get
\begin{equation*}\label{eq:painleve_x}x_{n+1}+x_n+x_{n-1}=\frac{2n+\alpha-\alpha(-1)^n}{8x_n}.\end{equation*}
Putting $y_n=a_n^2q^{1-n}$, we can factorize the equations: we
obtain
\begin{equation}\label{eq:painleve_y_evenodd}q^{n-\alpha}\left(y_ny_{n+1}+q^\alpha\right)\left(y_ny_{n-1}+q^\alpha\right)=\left\{\begin{array}{ll}q^\alpha-q^{-\alpha}y_n^2,\qquad&{\rm
for\: even}\:n\\1-y_n^2\qquad&{\rm for\:
odd}\:n.\end{array}\right.\end{equation} If we write
$u_n=q^{-\alpha}y_{2n}$ and $v_n=-y_{2n+1}$, substituting $n$ by
$2m$ resp. $2m+1$ in \reff{eq:painleve_y_evenodd}, we get
\begin{equation}\label{eq:uv_geen_c}\left\{\begin{array}{rcl}q^{2m}\left(1-u_mv_m\right)\left(1-u_mv_{m-1}\right)&=&1-u_m^2\\
q^{2m+1+\alpha}\left(1-u_mv_m\right)\left(1-u_{m+1}v_m\right)&=&1-v_m^2.\end{array}\right.\end{equation}
This set of equations could therefore be called a set of asymmetric
$q$-discrete Painlev\'{e} I equations ($\alpha$-q-P$_{\rm I}$).
However, in the next section we will obtain a more general form of
it.  Concerning the asymptotic behaviour of $a_n$ as $n$ tends to
infinity (with fixed $q$), it follows easily from
\reff{eq:recrel_a_n_evenodd} that $a_n$ tends to 0, and that
\begin{equation}\label{eq:a_n_asymp}\lim_{n\rightarrow\infty}\frac{a_{2n}^2}{q^{2n-1}}=\lim_{n\rightarrow\infty}y_{2n}=q^\alpha\qquad{\rm
and}\qquad\lim_{n\rightarrow\infty}\frac{a_{2n+1}^2}{q^{2n}}=\lim_{n\rightarrow\infty}y_{2n+1}=1,\end{equation}
or, equivalently, \[\lim_{n\rightarrow\infty}u_n=1\qquad{\rm
and}\qquad\lim_{n\rightarrow\infty}v_n=-1.\]The set of equations
\reff{eq:painleve_y_evenodd} has the singularity confinement
property: if $y_n=0$, then a singularity occurs for $y_{n+1}$ but
this singularity does not propagate infinitely: e.g. for $n$ even,
putting $y_n=\epsilon$ we obtain
\begin{eqnarray*}y_n&=&\varepsilon\\
y_{n+1}&=&-q^\alpha\left(1-q^{-n}\right)\frac{1}{\varepsilon}-y_{n-1}q^{-n}+\O(\varepsilon)\\
y_{n+2}&=&q^{\alpha-1}\left(1-q^{-n}\right)\frac{1}{\varepsilon}+\frac{y_{n-1}}{q}+\O(\varepsilon)\\
y_{n+3}&=&\frac{q^{n+1}-q^{-\alpha-2}}{1-q^n}\varepsilon+\O(\varepsilon^2)\\
y_{n+4}&=&q^{\alpha+2}\frac{1-q^n}{1-q^{n+\alpha+3}}y_{n-1}+\O(\varepsilon)
\end{eqnarray*} so the singularity is confined to
$y_{n+1},y_{n+2},y_{n+3}$; furthermore $y_{n+4}$ depends on the
value $y_{n-1}$ before the singularity. The same holds for odd $n$:
\begin{eqnarray*}y_n&=&\varepsilon\\
y_{n+1}&=&\left(q^{-n}-q^\alpha\right)\frac{1}{\varepsilon}-y_{n-1}q^{-n-\alpha}+\O(\varepsilon)\\
y_{n+2}&=&\left(q^{\alpha-1}-q^{-n-1}\right)\frac{1}{\varepsilon}+\frac{y_{n-1}}{q}+\O(\varepsilon)\\
y_{n+3}&=&q^\alpha\frac{q^{n}-q^{-3}}{1-q^{n+\alpha}}\varepsilon+\O(\varepsilon^2)\\
y_{n+4}&=&q^{2-\alpha}\frac{1-q^{n+\alpha}}{1-q^{n+3}}y_{n-1}+\O(\varepsilon).
\end{eqnarray*}
There are two (one for each parity) more critical cases that might
give rise to singularities, namely when $y_ny_{n-1}+q^\alpha=0$.
Running the singularity analysis gives, for even $n$,
\begin{eqnarray*}
y_n &=& -\frac{q^\alpha}{y_{n-1}}+\varepsilon\\
y_{n+1} &=&
\frac{q^{-n+\alpha}(y_{n-1}^2-1)}{y_{n-1}^2}\frac{1}{\varepsilon} +
\frac{q^{-n}+q^{-n}y_{n-1}^2-y_{n-1}^2}{y_{n-1}} +\O(\varepsilon)\\
\end{eqnarray*}
The coefficient of $1/\varepsilon$ however, is $\O(\varepsilon)$
itself as $y_{n-1}$ satisfies
$$
q^{n-\alpha-1}\left(y_ny_{n-1}+q^\alpha\right)\left(y_{n-2}y_{n-1}+q^\alpha\right)=1-y_{n-1}^2
$$
or, after substituing $y_n$,
$$
1-y_{n-1}^2 = \varepsilon
y_{n-1}q^{n-1-\alpha}(y_{n-1}y_{n-2}+q^\alpha).
$$
So there is no singularity to confine.

\begin{figure}
\centering
\includegraphics[angle=0,width=0.5\textwidth]{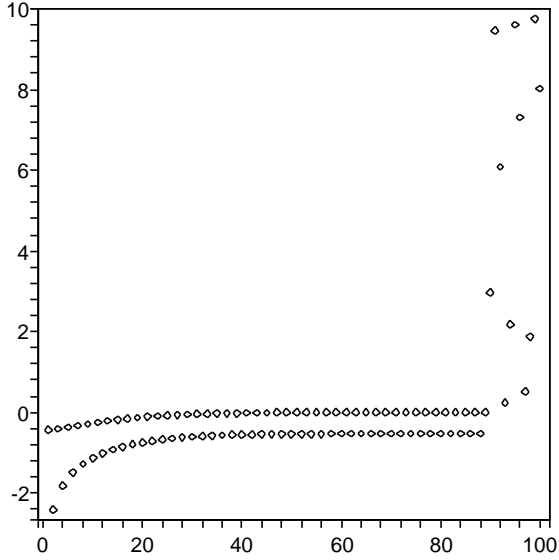}
\caption{The result of computing $\log|y_n|$ from
\reff{eq:painleve_y_evenodd} with $\alpha=5$, $q=0.9$, with an
accuracy of 200 digits. The odd-even dependence is clearly
visible.}\label{fig}
\end{figure}

\noindent One could hope to use the recurrence relations
(\ref{eq:painleve_y_evenodd}) to compute the recurrence
coefficients.  Given $y_0=0$, it would be nice to show that there is
a unique solution to \reff{eq:painleve_y_evenodd} which is positive
for all $n>0$.  The existence of such a solution is clear since
$y_n=a_n^2/q^{n-1}$ satisfies the given recurrence relations and is
obviously positive for $n>0$.  Its uniqueness however can only be
expected based on numerical experiments, but we have no proof of it
yet.  The method used in Section~\ref{section:stable} is probably
the key to proving this uniqueness.  Using \reff{eq:intro_an} and
the $q$-binomial theorem, it is clear that
\begin{equation*}y_1=a_1^2=\frac{\int_{-1}^1x^2\left(x^4q^4;q^4\right)_\infty|x|^\alpha d_qx}{\int_{-1}^1\left(x^4q^4;q^4\right)_\infty|x|^\alpha d_qx}=\frac{\left(q^{\alpha+1};q^4\right)_\infty}{\left(q^{\alpha+3};q^4\right)_\infty}.\end{equation*}
The knowledge of $y_0$ and $y_1$ and the recurrence relation
\reff{eq:painleve_y_evenodd} allow us to compute the $y_n$
recursively.  However, this turns out to be an unstable
 method. We computed $y_n$ to a certain accuracy for
some particular choices of $q$ and $\alpha$. In Figure~\ref{fig} we
plotted the values of $\log|y_n|$, calculated with an accuracy of
200 digits, with parameters $q=0.9,\alpha=5$. For this choice of
parameters, the computed values seem to satisfy the limit behaviour
as stated in \reff{eq:a_n_asymp} up to $n\sim 90$, but for larger
$n$ the sensibility of the non-linear equations on the initial
values destroys this behaviour. This is why
Section~\ref{section:stable} will be devoted to a more stable way of
computing the recurrence coefficients.

\section{Another modified discrete $q$-Freud
case}\label{section:c}

\begin{figure}
\centering
\includegraphics[angle=0,width=0.5\textwidth]{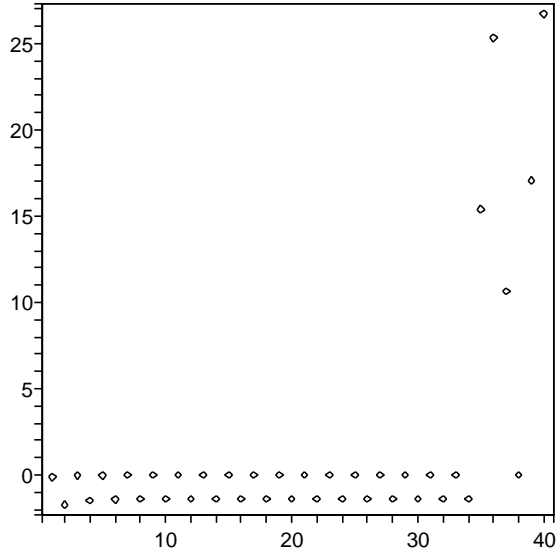}
\caption{The result of computing $\log|y_n|$ from
\reff{eq:c:painleve_y_evenodd} with $\alpha=2$, $q=0.5$, $c=-1/3$,
with an accuracy of 200 digits.}\label{fig2}
\end{figure}

In this section we consider a generalization of the previous case:
let the weight be given by
\begin{equation}\label{eq:weight_c_def}w(x)=\frac{|x|^\alpha\left(q^2x^2;q^2\right)_\infty\left(cq^2x^2;q^2\right)_\infty}{\left(1-q^4\right)^{\alpha/4}}\end{equation}
with $c\leq0$.  This weight satisfies the Pearson equation
\begin{equation*}\label{eq:c:pearson}w\left(\frac
xq\right)=\frac{\left(1-x^2\right)\left(1-cx^2\right)}{q^\alpha}w(x).\end{equation*}
If
\begin{equation}\label{eq:c:verband_a_c}c=-1+a\sqrt{1-q^4}\end{equation}
then
\begin{equation*}w\left(\sqrt[4]{1-q^4}x\right)\rightarrow|x|^\alpha
e^{-x^4-2ax^2}\qquad {\rm as }\: q\rightarrow 1,\end{equation*} so
this weight is a $q$-analogue of the Freud weight $|x|^\alpha
e^{-x^4-2ax^2}$. Obviously the choice $c=-1$ gives us the case
considered in the previous section. The choice $c=0$ leaves us with
a $q$-analogue of the weight $|x|^\alpha e^{-x^2}$, which was
already studied in \cite{walter} in the case $\alpha=0$.
\begin{lemma}\label{lemma:c:structure_relation}For even $n$, the orhogonal polynomials $p_n$ for the weight $w$ in \reff{eq:weight_c_def} satisfy the following structure
relation:
\begin{equation*}\label{eq:c:structure_even}\rd_qp_n(x)=\frac{\hat{B}_n}{1-q}p_{n-1}(x)+\frac{\hat{A}_n}{1-q}p_{n-3}(x)\end{equation*}
with
\begin{equation*}\label{eq:c:A_n_even}\hat{A}_n=-c\frac{a_na_{n-1}a_{n-2}}{q^{\alpha+n-3}}\end{equation*}
\begin{equation*}\label{eq:c:B_n_even}\hat{B}_n=\frac{a_n}{q^{\alpha+n-1}}\left(c+1-c\sum_{j=1}^{n+1}a_j^2+cq^2\sum_{j=1}^{n-2}a_j^2\right).\end{equation*}
For odd $n$, the polynomials $p_n$ satisfy
\begin{equation*}\label{eq:c:structure_odd}\rd_qp_n(x)=\frac{\hat{B}_n}{1-q}p_{n-1}(x)+\frac{\hat{A}_n}{1-q}p_{n-3}(x)+{\rm lower} \: {\rm order} \: {\rm terms}\end{equation*}
with
\begin{equation*}\label{eq:c:A_n_odd}\hat{A}_n=-c\frac{a_na_{n-1}a_{n-2}}{q^{\alpha+n-3}}-\left(1-q^{-\alpha}\right)\frac{a_{n-1}}{a_na_{n-2}}\end{equation*}
\begin{equation*}\label{eq:c:B_n_odd}\hat{B}_n=-c\frac{a_n}{q^{\alpha+n-1}}\left(\sum_{j=1}^{n+1}a_j^2-q^2\sum_{j=1}^{n-2}a_j^2\right)+\frac{1-q^{-\alpha}}{a_n}+(c+1)\frac{a_n}{q^{\alpha+n-1}}.\end{equation*}
\end{lemma}
\begin{proof}The proof follows the same steps as the proof of Lemma~\ref{lemma:structure_relation}.  The only difference is that in the Pearson equation, the quartic polynomial $1-x^4$ has now got to be replaced by $1-(c+1)x^2+cx^4$.\end{proof}
Performing the same manipulations as in the previous section leads
us to a recurrence relation satisfied by the $a_n$:
\begin{equation}\label{eq:c:recrel_a_n_evenodd}\left\{\begin{array}{rl}a_n^2\left(c+1-c\left(a_{n+1}^2+q^{-\alpha-n+1}a_n^2+q^2a_{n-1}^2-cq^{-2n-\alpha+3}a_{n+1}^2a_n^2a_{n-1}^2\right)\right)=\left(1-q^n\right)q^{\alpha+n-1},&n\:{\rm even}\\a_n^2\left(c+1-c\left(a_{n+1}^2+q^{-n+1}a_n^2+q^2a_{n-1}^2-cq^{-2n-\alpha+3}a_{n+1}^2a_n^2a_{n-1}^2\right)\right)=\left(q^{-\alpha}-q^n\right)q^{\alpha+n-1},&n\:{\rm odd.}\end{array}\right.\end{equation}
Putting again $y_n=a_n^2q^{1-n}$ we obtain a $q$-discrete
Painlev\'{e} I equation, of which \reff{eq:painleve_y_evenodd} is a
special case:
\begin{equation}\label{eq:c:painleve_y_evenodd}q^{n-\alpha}\left(-cy_ny_{n+1}+q^\alpha\right)\left(-cy_ny_{n-1}+q^\alpha\right)=\left\{\begin{array}{ll}\left(q^\alpha-y_n\right)\left(q^\alpha-cy_n\right)q^{-\alpha},\qquad&{\rm
for\: even}\:n\\\left(1-y_n\right)\left(1-cy_n\right),\qquad&{\rm
for\: odd}\:n.\end{array}\right.\end{equation} \noindent We can
rewrite these equations into a more familiar form for asymmetric
discrete Painlev\'e equations: denoting $u_n=y_{2n}q^{-\alpha}$ and
$v_n = cy_{2n+1}$ we get
\begin{equation}
\left\{ \begin{array}{rcl}
q^{2n}\left(1-u_nv_n\right)\left(1-u_nv_{n-1}\right) & =
& \left(1-u_n\right)\left(1-cu_n\right)\\
q^{2n+1+\alpha}\left(1-u_nv_n\right)\left(1-u_{n+1}v_n\right) &=&
\left(1-v_n\right)\left(1-v_n/c\right).
\end{array}\right.\label{eq:asymm_vorm}
\end{equation}  Besides $q$, this set of equations contains the
parameters $\alpha$ and $c$, and it is the most general form of the
asymmetric $q$-discrete Painlev\'e I equations we obtain. Hence we
call this set $\alpha$-q-P$_{\rm I}$.  It is interesting to notice
that these modifications w.r.t. the previous section 'survive' in
the limit case $q\rightarrow1$: if we put
$x_n=\frac{a_n^2}{\sqrt{1-q^4}}$ and we let $q$ tend to 1, we obtain
\begin{equation*}\label{eq:c:painleve_x}x_{n+1}+x_n+x_{n-1}=\frac{2n+\alpha-\alpha(-1)^n}{8x_n}-a,\end{equation*}
which is still an instance of the discrete Painlev\'{e} I equation
\reff{eq:dPI}.  Here $a$ and $c$ are related as in
\reff{eq:c:verband_a_c}. The introduction of the additional
parameter $c$ has no influence on the limit behaviour of the
recurrence coefficients $a_n$: \reff{eq:a_n_asymp} still holds. To
show this we distinguish the following cases.  We only mention the
proof for even $n$, the odd $n$ case being analogous.
\begin{itemize}
\item $c=-1$ is the case from the previous section.
\item $0>c>-1$: for even $n$, \reff{eq:c:recrel_a_n_evenodd} gives
\[\frac{a_n^4}{q^{2n-2}}\leq\frac{q^{2\alpha}(1-q^n)}{-c}\]
so $a_n\rightarrow 0$ and $y_n=\frac{a_n^2}{q^{n-1}}$ is bounded.
Denoting $A=\limsup y_n$, \reff{eq:c:recrel_a_n_evenodd} gives the
quadratic equation $A\left(c+1-\frac{cA}{q^\alpha}\right)=q^\alpha$.
This equation has a negative solution $A=q^\alpha/c$ which we can
reject, and a positive solution $A=q^\alpha$.  The same argument can
be used on $B=\liminf y_n$, and one obtains \reff{eq:a_n_asymp}.
\item $c<-1$, $n$ even: from \reff{eq:c:recrel_a_n_evenodd} we get
\[\frac{a_n^2}{q^{\alpha+n-1}}\left(c+1-c\frac{a_n^2}{q^{\alpha+n-1}}\right)\leq 1-q^n<1.
\]
Looking at this as a quadratic inequality for
$y_n=\frac{a_n^2}{q^{n-1}}$ we obtain
\[\frac{y_n}{q^\alpha}\left(c+1-\frac{cy_n}{q^\alpha}\right)<1\] and hence $y_n\in\left(1/c,1\right)$.  A similar argument with $\limsup$ and $\liminf$ yields \reff{eq:a_n_asymp}.
\item $c=0$: then the result follows immediately from
\reff{eq:c:recrel_a_n_evenodd}.
\end{itemize} Once more, the singularity
confinement property is fulfilled: taking $y_n=\epsilon$ for an even
$n$, one finds that $y_{n+1}$ and $y_{n+2}$ are $\O
(1/\varepsilon)$, $y_{n+3}$ is $\O (\varepsilon)$ and
\[y_{n+4}=\frac{q^{2+\alpha}}{c\left(1-q^{n+\alpha+3}\right)}\left[cy_{n-1}(1-q^n)-(1+c)(1-q^{-2})\right]+\O(\epsilon).\]
Now it is not obvious that the constant term is nonzero: if
$y_{n-1}=\frac{(1+c)(1-q^{-2})}{c(1-q^n)}$ then this is not yet a
proof that the singularity confinement property holds. For odd $n$ a
similar argument yields
\[y_{n+4}=\frac{q^{2-\alpha}}{c\left(1-q^{n+3}\right)}\left[cy_{n-1}(1-q^{n+\alpha})-(1+c)(1-q^{-2})q^{\alpha}\right]+\O(\epsilon)\]
and a problem could arise if
$y_{n-1}=\frac{q^\alpha(1+c)(1-q^{-2})}{c(1-q^{n+\alpha})}$.
 However (consider the case $n$ is even, the odd $n$ case being completely analogous) the zero at $y_{n+4}$ gives rise to singularities in
 $y_{n+5},y_{n+6},y_{n+7}$, but the singularity vanishes in
 $y_{n+8}$.  Moreover, for the same reasons as above, this $y_{n+8}$
 can only be zero if
 $y_{n+3}=\frac{1+c}{c}\frac{1-q^{-2}}{1-q^{n+4}}$ (which is nonzero since $|q|<1$ and $c\neq-1$, the case $c=-1$ being considered in the previous section), while the
 computation gave $y_{n+3}=\O(\epsilon)$.  So this assures
 us that even in the worst case, the singularity is confined to
 $y_{n+1},\ldots,y_{n+7}$.
Again, as in Section~\ref{section:nog_geen_c}, there are no
singularities arising from $-cy_ny_{n-1}+q^\alpha=0$ due to fine
cancellations.
 \\
\noindent Concerning the use of the recurrence relations to compute
the recurrence coefficients $a_n$ of the orthogonal polynomials, we
now start with $y_0=0$ and
$$y_1=a_1^2=\frac{\sum_{k=0}^\infty \frac{q^{k(\alpha+3)}}{(q^2;q^2)_k(cq^2,q^2)_k}}{\sum_{k=0}^\infty \frac{q^{k(\alpha+1)}}{(q^2;q^2)_k(cq^2,q^2)_k}}.$$ The same remark as in the previous
section holds: this method is unstable, as is shown in
Figure~\ref{fig2}. We refer to Section~\ref{section:stable} for a
stable computation method.
\begin{remark}We restricted ourselves to the choice $c<0$, with the case $c=0$ being a limiting case where the Painlev\'e equations give explicit expressions for e.g. $y_n$.  However, numerical experiments strongly suggest that also for $0<c<1$ (in which case the weight $w(x/q)$ is still positive on $[-1,1]$), the results on the limit behaviour \reff{eq:a_n_asymp} still hold, and the singularity confinement property is still valid as well.  Even for $c\geq 1$, where for small $n$ we observe that $y_n$ can be negative, we still observe that $y_n$ is bounded and oscillating, however the even and odd subsequences do not tend to $q^\alpha$ and 1 respectively anymore, but to $q^\alpha/c$ and $1/c$, respectively.\end{remark}

\section{Connection to $\alpha$-q-P$_{{\rm
V}}$}\label{section:painleve5} It is worth noticing that the
obtained recurrence relations can be seen as a limiting case of the
set of asymmetric discrete Painlev\'e equations $\alpha$-q-P$_{{\rm
V}}$.  This set can be found as the discrete Painlev\'e equation
connected to the affine Weyl group $E_6^q$ in Sakai's
classification(see, e.g., \cite{gramram}), or as $\alpha$-q-P$_{{\rm
V}}$ in \cite{walter}. The equations are
\[\left(1-u_nv_n\right)\left(1-u_nv_{n-1}\right)=\frac{\left(u_n-1/p\right)\left(u_n-1/r\right)\left(u_n-1/s\right)\left(u_n-1/t\right)}{\left(u_n-b\rho_n\right)\left(u_n-\rho_n/b\right)}\]
\[\left(1-u_nv_n\right)\left(1-u_{n+1}v_n\right)=\frac{\left(v_n-p\right)\left(v_n-r\right)\left(v_n-s\right)\left(v_n-t\right)}{\left(v_n-aw_n\right)\left(v_n-w_n/a\right)}\]
with $prst=1$.  Consider this system with the particular choice of
parameters
\[p=1,r=c,s=\kappa,t=\frac{1}{c\kappa},b=c\kappa,\rho_n=q^{2n},a=\kappa,w_n=q^{2n+\alpha+1}.\]
Then, letting $\kappa$ tend to 0, we obtain as a limit the set of
equations \reff{eq:asymm_vorm}.

\section{A stable method for computing the recurrence coefficients}\label{section:stable}
\begin{figure}
\centering
\includegraphics[angle=0,width=0.5\textwidth]{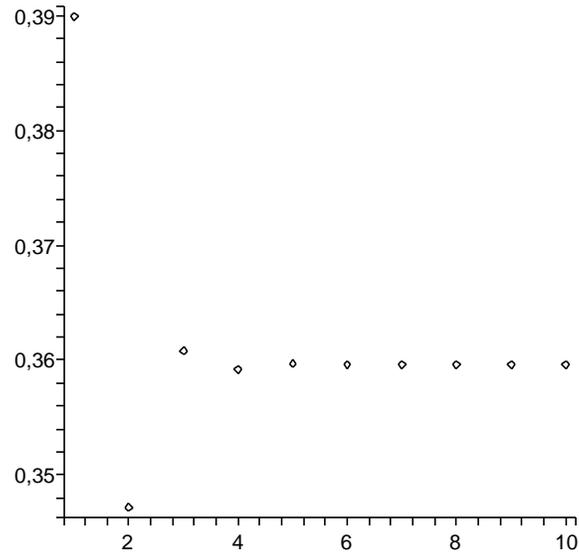}
\caption{Approximation of the recurrence coefficient $y_1$ for
$\alpha=2,q=9/10,c=-1/2$; on the horizontal axis is the number of
times we applied the operator $T$ to the starting value
$(0,0)$.}\label{fig:stabiel1}
\end{figure}

\begin{figure}
\centering
\includegraphics[angle=0,width=0.5\textwidth]{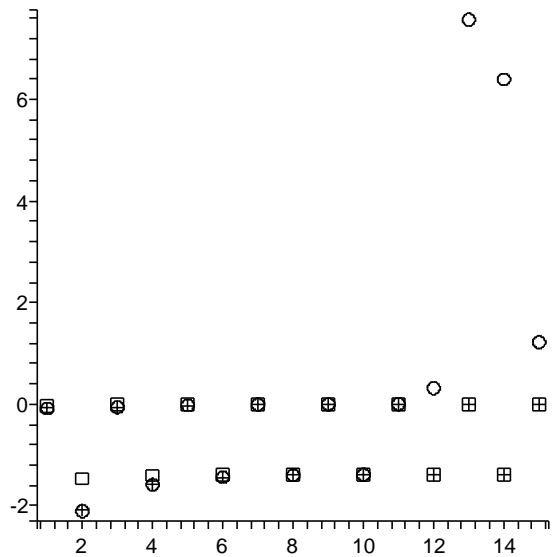}
\caption{$\log|y_n|$ with $\alpha=2,q=1/2,c=-5/2$, computed using
the recurrence relations \reff{eq:c:painleve_y_evenodd} with 20
digits accuracy (circles), and the approximation obtained after
applying once (squares), or three times (crosses) the operator $T$
of Section~\ref{section:stable} to the starting value
$(0,0)$.}\label{fig:allebei}
\end{figure}

\noindent A more stable way for computing the recurrence values is
by writing \reff{eq:c:painleve_y_evenodd} as a system of quadratic
equations in $\tilde{u}_n:=y_{2n}$ and $\tilde{v}_n:=y_{2n+1}$, in
the same way as we obtained \reff{eq:asymm_vorm}:
$$
\left\{
\begin{array}{rcl}\tilde{u}_n^2q^{-\alpha}\left(c^2q^{2n}\tilde{v}_n\tilde{v}_{n-1}-c\right)+\tilde{u}_n\left(-cq^{2n}(\tilde{v}_n+\tilde{v}_{n-1})+c+1\right)+q^{\alpha}\left(q^{2n}-1\right) &=&0\\
\tilde{v}_n^2\left(c^2q^{2n+1-\alpha}\tilde{u}_n\tilde{u}_{n+1}-c\right)+\tilde{v}_n\left(-cq^{2n+1}(\tilde{u}_n+\tilde{u}_{n+1})+c+1\right)+\left(q^{2n+1+\alpha}-1\right)
&=&0.\end{array}\right.
$$
Computing the discriminants and opting for their positive roots (as
we know, the $\tilde{u}_n$ and $\tilde{v}_n$ are related to the
recurrence coefficients which are positive) in the expression of
$\tilde{u}_n$ and $\tilde{v}_n$, we find
$$
\left\{ \begin{array}{rcl} \tilde{u}_n &=&
f_n(-cq^{2n}(\tilde{v}_n+\tilde{v}_{n-1})+c+1,c^2q^{2n}\tilde{v}_n\tilde{v}_{n-1}-c)\\
\tilde{v}_n & = &
g_n(-cq^{2n+1}(\tilde{u}_n+\tilde{u}_{n+1})+c+1,c^2q^{2n+1-\alpha}\tilde{u}_n\tilde{u}_{n+1}-c)
\end{array}\right.
$$
with
$$
\left\{ \begin{array}{rcl} f_n(x,y) &=&
\frac{-x+\sqrt{x^2+4(1-q^{2n})y}}{2q^{-\alpha}y}\\
g_n(x,y)&=& \frac{-x+\sqrt{x^2+4(1-q^{2n+1+\alpha})y}}{2y}.
\end{array}\right.
$$
We now define an operator $T$, acting on the space of double
non-negative rows:
$$
T\left(\begin{array}{rcl} \xi&=& (\xi_0, \xi_1, \ldots)\\ \eta & =
&(\eta_0, \eta_1, \ldots) \end{array} \right) = \left(
\begin{array}{c}\xi'\\ \eta' \end{array} \right)
$$
where $$\xi'_0 = 0,\,\,\, \xi'_n =
f_n(-cq^{2n}(\eta_n+\eta_{n-1})+c+1,c^2q^{2n}\eta_n\eta_{n-1}-c)
\textrm{ for } n>0$$ and $$\eta'_n
=g_n(-cq^{2n+1}(\xi_n+\xi_{n+1})+c+1,c^2q^{2n+1-\alpha}\xi_n\xi_{n+1}-c)
\textrm{ for } n \geq 0.$$ The solution of the equation which comes
from our recurrence coefficients, denoted by
$(\tilde{u},\tilde{v})$, now coincides with a fixed point of this
operator $T$. We know that $(\tilde{u},\tilde{v}) \geq (0,0)$
(inequalities should be interpreted termwise: as holding between any
two elements on corresponding positions). Unleashing the operator
$T$ (on any inequality we will use) reverses inequalities.
\begin{lemma}
If $(\xi,\eta) \leq (a,b)$, then $ T(\xi,\eta) \geq T(a,b)$.
\end{lemma}
\begin{proof}
It is not hard to check that all partial derivatives of $f_n$ and
$g_n$ are negative in the region where we need them.  Which region
this is, depends on the parameter $c$.
\end{proof}
We now construct the following sequence by repeatedly applying $T$
to $(0,0)$. This sequence has increasing subsequence $T^{2k}(0,0)$
and decreasing subsequence $T^{2k+1}(0,0)$. We denote their limits
by $(\xi^-,\eta^-)$ and $(\xi^+,\eta^+)$, respectively. By
continuity of $T$ we get that $T(\xi^-,\eta^-)=(\xi^+,\eta^+)$ and
$T(\xi^+,\eta^+)=(\xi^-,\eta^-)$. Furthermore, any fixed point
$(\xi^*,\eta^*)$ with positive components has to obey
$(\xi^-,\eta^-)\leq (\xi^*,\eta^*) \leq (\xi^+,\eta^+)$. If one is
able to show that $(\xi^-,\eta^-)= (\xi^+,\eta^+)$, one proves there
is only one fixed point, which necessarily coincides with the
solution of the discrete Painlev\'e equation coming from the
recurrence coefficients. We did not find this proof (yet), but
numerical experiments strongly suggest that indeed $(\xi^-,\eta^-)=
(\xi^+,\eta^+)$. This would also be enough to prove the uniqueness
of positive solutions to \reff{eq:c:painleve_y_evenodd} or
\reff{eq:painleve_y_evenodd} with $y_0=0$, as mentioned in
Section~\ref{section:nog_geen_c}. See Figure~\ref{fig:stabiel1} for
the result using a particular choice of parameters, and
Figure~\ref{fig:allebei} for a comparison between the two methods.

\begin{remark}If $c=0$ then
\reff{eq:c:painleve_y_evenodd} immediately gives the exact
recurrence coefficients: $y_n=q^\alpha-q^{n+2\alpha}$ for even $n$,
and $y_n=1-q^{n+\alpha}$ for odd $n$, which gives
$$a^2_n = \left\{ \begin{array}{ll} q^{n+\alpha-1} (1-q^{n+\alpha}) & \textrm{for even }n\\ q^{n-1}(1-q^{n+\alpha}) & \textrm{for odd }n. \end{array}\right.$$
\end{remark}

\section{Conclusions}\label{section:conclusions}
\noindent We used a generalized $q$-Freud weight in order to find a
$q$-discrete Painlev\'e equation. What we found was an asymmetric
form of d-P$_{\rm I}$ which has (as we believe) never been described
before. The asymmetric form emerges in a very natural way from the
parity of the orthogonal polynomials associated to the weight. We
showed its relation to $\alpha-q$P$_{\rm V}$ and gave a stable
method to compute the recurrence coefficients for the orthonormal
polynomials associated to this weight.

\noindent The technique used in this paper is not sufficient to
handle weights that are not even. In these cases, one has a non-zero
$b_n$ entering the system, which could lead to new results. This is
a possible interest for future research.

\section{Acknowledgements}\label{section:ackn}
This research was supported by Belgian Interuniversity Attraction
Poles Programme P6/02, K.U.Leuven projects OT/04/21 and OT/08/33.

\end{document}